\documentclass[12pt]{amsart}
\usepackage{latexsym}
\usepackage{amssymb, amsmath,mathrsfs}
\usepackage[ansinew]{inputenc}
\usepackage[T1]{fontenc}
\usepackage{amssymb}
\usepackage{amsthm}
\usepackage{amsfonts}
\usepackage{amsmath}
\usepackage{lmodern}
\usepackage[margin=1in]{geometry}
\usepackage{comment}
\usepackage{color}

\usepackage{titlesec}
\titlelabel{\thetitle.\,\:}
\titleformat*{\section}{\bf\center}

\newcommand{\R}{\mathbb{R}}

\newcommand{\supp}{\textrm{supp}\,}
\newcommand{\ch}{\textrm{ch}\,}

\theoremstyle{definition}
\newtheorem{rema}{Remark}
\newtheorem{defi}{Definition}
\theoremstyle{plain}
\newtheorem{lem}{Lemma}
\newtheorem{theo}[lem]{Theorem}

\theoremstyle{definition}

\theoremstyle{plain}

\newenvironment{proofof}[1]
{\noindent{\it Proof of #1. }}{\hfill $\Box$\par\vspace{2.5mm}}

\begin{document}

\title{\bf Two-Weight $Tb$ Theorems for Well-Localized Operators}
\date{}
\author[K. Bickel]{Kelly Bickel$^1$}
\address{Kelly Bickel, Department of Mathematics\\
Bucknell University\\
701 Moore Ave\
Lewisburg, PA 17837}
\email{kelly.bickel@bucknell.edu}
\thanks{$1.$ Research supported in part by National Science Foundation
DMS grant \#1448846.}

\author[T. Korhonen]{Taneli Korhonen$^2$}
\address{Taneli Korhonen, Department of Mathematics\\
University of Eastern Finland, Department of Physics and Mathematics,\\
P.O.~Box 111, 80101 Joensuu, Finland}
\email{taneli.korhonen@uef.fi}
\thanks{$2.$ Research supported in part by Faculty of Science and Forestry of University of Eastern Finland.}

\author[B. D. Wick]{Brett D. Wick$^3$}
\address{Brett D. Wick, Department of Mathematics\\
Washington University in St. Louis\\
One Brookings Drive\\
 St. Louis, MO 63130-4899 }
\email{wick@math.wustl.edu}
\thanks{$3.$ Research supported in part by National Science Foundation
DMS grant \#1500509 and \# 1800057.}

\maketitle

\renewcommand{\thefootnote}{}


\begin{abstract}
This paper first defines operators that are ``well-localized'' with respect to  a pair of accretive functions and establishes a global two-weight $Tb$ theorem for such operators.
Then it defines  operators that are ``well-localized'' with respect to a pair of accretive systems and establishes a local two-weight $Tb$ theorem for them.  The proofs combine recent $Tb$ proof techniques with arguments used to prove earlier $T1$ theorems for well-localized operators.
\end{abstract}

\section{Introduction}
Over the past several decades, researchers have proved a number of important $Tb$ theorems showing that the boundedness of  Calder\'on-Zygmund operators can be deduced from testing on certain functions $b$.
David,  Journ\'e, and Semmes proved the first global $Tb$ theorem in 1985  \cite{DJS1985}; they showed that for sufficiently nice (accretive) functions $b$ and $c$, a Calder\'on-Zygmund operator $T:L^2(m)\rightarrow L^2(m)$ is bounded precisely when $M_b TM_c$ is weakly bounded and $Tb, T^*c \in \text{BMO}.$ This result was generalized to nonhomogeneous settings in both \cite{HM2012B, NTV2003}.

Meanwhile in 1990, Christ established a local $Tb$ theorem in the homogeneous setting by showing  that $T:L^2(\mu)\rightarrow L^2(\mu)$ is bounded if $\| Tb_Q\|_{L^\infty(\mu)}$ and $\| b_{Q}\|_{L^\infty(\mu)}$ are uniformly bounded for a system of accretive functions $\{b_Q\}$ \cite{Christ1990}. This theorem was generalized to the nonhomogeneous settings in \cite{HM2012, NTV2002}. In the homogeneous setting, alternate $L^p$--rather than $L^{\infty}$--testing conditions have also been studied extensively, see \cite{AHMTT02, AR13, AY09,  Hof07, HN12, LV12, TY09}, and some $L^p$ testing conditions have even been extended to the nonhomogeneous setting \cite{HM2012, LM2016}.

This paper connects this rich field of $Tb$ theorems to the setting of well-localized operators, which were studied in \cite{BW16, BCTW18,NTV2008}. Well-localized operators are closely connected to band, or almost-diagonal, operators. Indeed, in both \cite{ BCTW18, NTV2008}, the authors showed that the boundedness of band operators, such as Haar shifts of a fixed complexity, is equivalent to the boundedness of certain well-localized operators. Motivated by such connections, the authors in  \cite{BW16, BCTW18, NTV2008} established various $T1$ theorems for well-localized operators.

In this paper, we  extend the results from \cite{NTV2008} by establishing both global and local $Tb$ theorems for associated ``well-localized'' operators. In both settings, we let $\mu, \nu$ denote Borel measures that are nonnegative and finite on dyadic cubes $Q \in \mathcal{D}.$  Then in the global setting, we consider pairs of functions $b_1 \in L^{\infty}(\mu), b_2 \in L^{\infty}(\nu)$ with averages $|\langle b_1\rangle^{\mu}_Q |, |\langle b_2\rangle^{\nu}_Q |  \gtrsim 1$ for all dyadic cubes $Q\in \mathcal{D}$. Such $\{b_1,b_2\}$ are called $(\mu, \nu)$-weakly accretive. In Definition \ref{def:WLO}, we explain what it means for an operator $T$ to be well-localized (with radius $r$) with respect to such  $\{b_1, b_2\}$, and then we establish the following theorem:

\begin{theo} \label{thm:GTB} If $T$ is a $\{b_1, b_2\}$-well-localized operator with radius $r$ satisfying
\begin{itemize}
\item[(a)]  $\| 1_Q T(b_1 1_Q)\|^2_{L^2(\nu)} \lesssim \mu(Q)$ and $\| 1_Q T^*(b_2 1_Q)\|^2_{L^2(\mu)} \lesssim \nu(Q)$ for all $Q \in \mathcal{D},$
\item[(b)] $ |\langle T(1_Q b_1), 1_R b_2 \rangle_{\nu}| \lesssim \| 1_Q b_1 \|_{L^2(\mu)} \| 1_R b_2 \|_{L^2(\nu)}$ for all $Q,R \in \mathcal{D}$ satisfying $ 2^{-r} \ell(Q) \le \ell(R) \le 2^r \ell(Q)$,
\end{itemize}
then $T: L^2(\mu) \rightarrow L^2(\nu)$ is bounded. \end{theo}

This theorem is very much in the flavor of the $T1$ theorems from \cite{BW16, BCTW18, NTV2008}, and the proof adapts both $Tb$ arguments from \cite{NTV2003} and well-localized arguments from \cite{NTV2008}. For a complete explanation of the notation and further details, see Sections \ref{sec:GTB} and \ref{sec:GTBP}.

In the local setting, we prove a similar theorem, but with testing on accretive systems $\{b_Q\}_{Q \in \mathcal{D}}$, indexed by the dyadic lattice $\mathcal{D}$. This situation is more complicated and we adapt local $Tb$ arguments from both \cite{NTV2002} and \cite{HM2012}. Our proof techniques require additional assumptions on the accretive systems and their relationships to the measures as well as an additional testing condition that is trivial when the measures are doubling. The definition of a well-localized operator also requires a restrictive extra condition given in \eqref{eqn:csc}. With those assumptions, we prove Theorem \ref{theo: Tb-Theorem}, a local $Tb$ theorem that is similar  to Theorem \ref{thm:GTB} given above. The details can be found in Sections \ref{sec:LTB} and \ref{sec:LTBP}.

\section{Global $Tb$ Theorem} \label{sec:GTB}
Let $\mathcal{D}$ be the standard dyadic lattice in $\R^n$. In what follows, for each cube $Q\in \mathcal{D}$, $\ell(Q)$ denotes the side length of $Q$ and $\ch Q$ denotes the set of children of $Q$, namely the set of cubes $Q' \in \mathcal{D}$ satisfying $Q' \subset Q$ and $\ell(Q') = \ell(Q)/2.$ Similarly, $\ch^r Q$ denotes the set of cubes $Q' \subseteq Q$ with $\ell(Q') = 2^{-r} \ell(Q).$
Furthermore, $Q^{(1)}$ denotes the parent of $Q$ and $Q^{(r)}$ denotes the ancestor of $Q$ of order $r$, namely  $Q^{(r)}$ is the unique cube satisfying $Q \subseteq Q^{(r)}$ and $\ell(Q^{(r)})= 2^r \ell(Q).$ For a Borel measure $\mu$ and $f \in L^2(\mu)$, denote the average of $f$ over a cube $Q\in \mathcal{D}$ by
	\[
	\langle f\rangle^\mu_Q = \mu(Q)^{-1} \int_Q f\,d\mu.
	\]
To avoid dividing by zero, if $\mu(Q) =0$, set $\langle f\rangle^\mu_Q  \equiv 0.$ However, in the later proofs and formulas, for simplicity we will make the standard assumption that $\mu(Q) \ne 0 $ for all $Q \in \mathcal{D}.$ Given Borel measures $\mu, \nu$ on $\mathbb{R}^n$, we can define the testing functions.

\begin{defi}  \label{def:acc} A function $b$ is \emph{$\mu$-weakly accretive} if $b \in L^{\infty}(\mu)$
and the weighted averages of $b$ satisfy  $|\langle b\rangle^{\mu}_Q |  \gtrsim 1$
for all $Q\in \mathcal{D}$, with implied constant independent of $Q$. If $b_1$ is $\mu$-weakly accretive and $b_2$ is $\nu$-weakly accretive, then the pair $\{b_1, b_2\}$ is called \emph{$(\mu, \nu)$-weakly accretive.}
\end{defi}

Given a $\mu$-weakly accretive $b$, one can define the following expectations and martingale differences for each $f \in L^2(\mu)$ and $Q \in \mathcal{D}$:
\[
E^{b}_Q f := \frac{\langle f \rangle^\mu_Q}{ \langle b \rangle^\mu_Q} 1_Q b \ \ \ \text{ and } \ \ \ \Delta^b_Q f := \sum_{Q' \in \ch Q} E^b_{Q'} f - E^b_Q f.
\]
In what follows, any function in the range space $\Delta^b_Q L^2(\mu)$ is called a $\mu$-Haar function associated to $Q$. These functions are orthogonal to constants and are supported on $Q$. An arbitrary function in $\Delta^b_Q L^2(\mu)$ will be written as $h^b_Q.$   Furthermore, the operators $E^{b}_Q $ and $\Delta^b_Q$ are projections and give useful decompositions of $L^2(\mu)$ functions.

\begin{lem}\cite[pp. 192-193]{NTV2003} \label{lem:decomp0} Let $b$ be a $\mu$-weakly accretive function and let $f \in L^2(\mu).$ Then for each $d\in \mathbb{Z}$,
\[ f = \sum_{\substack{ Q \in \mathcal{D} \\ \ell(Q) \le 2^d}} \Delta^b_Q f + \sum_{\substack{Q \in \mathcal{D} \\ \ell(Q) = 2^d}} E^b_Q f\]
with convergence in $L^2(\mu)$. Moreover, the following estimates hold
\begin{itemize}
\item[(i)] $\displaystyle \sum_{\substack{ Q \in \mathcal{D} \\ \ell(Q) \le 2^d}} \|  \Delta^b_Q f \|^2_{L^2(\mu)}  + \sum_{\substack{Q \in \mathcal{D} \\ \ell(Q) = 2^d}} \| E^b_Q f\|^2_{L^2(\mu)}
\lesssim \|f \|^2_{L^2(\mu)}$,
\item [(ii)]
$\displaystyle \sum_{Q \in \mathcal{D}}  \|  (\Delta^b_Q)^*  f \|^2_{L^2(\mu)} \lesssim \|f \|^2_{L^2(\mu)}.$
\end{itemize}
\end{lem}
A basic estimate using the properties of $b$ shows that each $\| E^b_{Q} f \|_{L^2(\mu)} \lesssim \| f 1_{Q} \|_{L^2(\mu)}$, and thus
\begin{equation} \label{eqn:bd1} \Big \|  \sum_{\substack{ Q \in \mathcal{D} \\ \ell(Q) \le 2^d}}  \Delta^b_Q f \Big \|_{L^2(\mu)} \lesssim \|f\|_{L^2(\mu)}.\end{equation}
Given this setup, we can define the well-localized operators. Specifically, for a pair of $(\mu, \nu)$-weakly accretive functions $\{b_1, b_2\}$, we say that an operator $T$ acts formally from $L^2(\mu)$ to $L^2(\nu)$ with respect to $\{b_1, b_2\}$ if the bilinear form
\[ \langle T(b_1 1_Q), b_2 1_R \rangle_{L^2(\nu)}\]
is well defined for all $Q,R \in \mathcal{D}$.
\begin{defi} \label{def:WLO} Let $T$ be an operator acting formally from $L^2(\mu)$ to $L^2(\nu)$ with respect to $\{b_1, b_2\}$. Then $T$ is \emph{lower triangularly localized with respect to $\{b_1, b_2\}$ with radius $r$} if there exists an $r\ge 0$ such that for all $Q,R \in \mathcal{D}$ with $\ell(R) \le 2 \ell(Q)$,
\[ \langle T(b_1 1_Q), h^{b_2}_R \rangle_{\nu} =0 \qquad \forall \ h^{b_2}_R \in \Delta^{b_2}_R L^2(\nu),\]
whenever $R \not \subseteq Q^{(r)}$ or if $ \ell(R) \le 2^{-r} \ell (Q)$ and $R\not \subseteq Q.$  We say that $T$ is \emph{$\{b_1, b_2\}$-well-localized with radius $r$} if both $T$ and $T^*$ are lower triangularly localized with respect to $\{b_1, b_2\}$ with radius $r$. For $T^*$, the roles of $\mu$ and $\nu$ and $b_1$ and $b_2$ are switched.
\end{defi}
Then, as mentioned in the introduction, the following theorem can be proved in a way similar to the standard situation discussed in \cite[Theorem 2.3]{NTV2008}:\\

\noindent \textbf{Theorem \ref{thm:GTB}.} \emph{Let $\mu, \nu$ be Borel measures on $\mathbb{R}^n$ and let $\{ b_1, b_2\}$ be $(\mu,\nu)$-weakly accretive. Let $T$ be a $\{b_1, b_2\}$-well-localized operator with radius $r$ satisfying
\begin{itemize}
\item[(a)]  $\| 1_Q T(b_1 1_Q)\|^2_{L^2(\nu)} \lesssim \mu(Q)$ and $\| 1_Q T^*(b_2 1_Q)\|^2_{L^2(\mu)} \lesssim \nu(Q)$ for all $Q \in \mathcal{D},$
\item[(b)] $ |\langle T(1_Q b_1), 1_R b_2 \rangle_{\nu}| \lesssim \| 1_Q b_1 \|_{L^2(\mu)} \| 1_R b_2 \|_{L^2(\nu)}$ for all $Q,R \in \mathcal{D}$ satisfying $ 2^{-r} \ell(Q) \le \ell(R) \le 2^r \ell(Q)$.
\end{itemize}
Then $T: L^2(\mu) \rightarrow L^2(\nu)$ is bounded.} \\

\begin{rema}
\label{rem:MD} Testing conditions of the form $\| 1_Q T(b_1 1_Q)\|^2_{L^2(\nu)} \lesssim \| b_1 1_Q\|^2_{L^2(\mu)}$ can be used instead of (a). In particular, since $b_1 \in L^{\infty}(\mu)$, we would immediately obtain testing condition (a) by
\[ \| 1_Q T(b_1 1_Q)\|^2_{L^2(\nu)} \lesssim \| b_1 1_Q\|^2_{L^2(\mu)} \lesssim \mu(Q).\]
Similarly, a simple argument using  testing condition (b) and the definition of our martingale differences
shows that
\[ |\langle T (\Delta^{b_1}_Q f), \Delta^{b_2}_R g \rangle_{\nu}| \lesssim \| \Delta^{b_1}_Q f \|_{L^2(\mu)} \| \Delta^{b_2}_R g \|_{L^2(\nu)}\]
for all $f \in L^2(\mu)$ and $g\in L^2(\nu)$ and all $Q,R \in \mathcal{D}$ satisfying $ 2^{-r} \ell(Q) \le \ell(R) \le 2^r \ell(Q)$.

\end{rema}

\section{Proof of Theorem \ref{thm:GTB}} \label{sec:GTBP}

The proof uses the following well-known theorem.

\begin{theo}[Dyadic Carleson embedding theorem]\label{theo: Carleson embedding}
If $\mu$ is a Borel measure and if $(a_Q)_{Q \in \mathcal{D}}$ is a $\mu$-Carleson sequence, i.e. if
\[ \sum_{Q\subseteq R} a_Q \lesssim \mu(R) \quad \forall R \in \mathcal{D},
  \ \text{ then } \
    \sum_{Q \in \mathcal{D}} a_Q \left|\langle f\rangle^\mu_Q\right|^2 \lesssim \|f\|^2_{L^2(\mu)}
    \quad \forall f \in L^2(\mu).
    \]
\end{theo}

\noindent \begin{proofof}{Theorem~\ref{thm:GTB}} Fix $f \in L^2(\mu)$ and $g \in L^2(\nu)$ and without loss of generality, assume that they are compactly supported. Then, there is an integer $d$ and cubes $Q_1, \dots, Q_{2^n} \in \mathcal{D}$ with no common ancestors such that $\ell(Q_j) = 2^d$ and $\supp f, \supp g \subseteq \cup Q_j.$ By Lemma \ref{lem:decomp0}, we can write
\[
\begin{aligned}
f & = \sum_{Q \subseteq \cup Q_j} \Delta^{b_1}_Q f + \sum_{j=1}^{2^n} E^{b_1}_{Q_j} f = f_1 + f_2 \\
g &= \sum_{R \subseteq \cup Q_k} \Delta^{b_2}_R g + \sum_{k=1}^{2^n} E^{b_2}_{Q_k}  g = g_1 + g_2.
\end{aligned}
\]
By duality, it suffices to show that $|\langle Tf,g\rangle_{\nu}| \le C \| f\|_{L^2(\mu)} \| g \|_{L^2(\nu)}.$ We break the inner product into the following four terms
\[ S_1 + S_2 +S_3 +S_4 = \langle Tf_1,g_1\rangle_{\nu} + \langle Tf_2,g_1\rangle_{\nu} + \langle Tf_1,g_2\rangle_{\nu} + \langle Tf_2,g_2\rangle_{\nu}\]
and handle them separately. We leave $S_1$ for later. First consider $S_2$ and observe that if $j \ne k$ and  $R \subseteq Q_k$, then $\ell(R) \le \ell(Q_j)$ and $R \not \subseteq Q^{(r)}_j$. Then the definition of well-localized implies that
\[ \langle T(E^{b_1}_{Q_j} f),  \Delta^{b_2}_R g  \rangle_{\nu} = 0.\]
This means that we can control $S_2$ by
\[\begin{aligned}
|S_2| &\le \sum_{j=1}^{2^n} \left | \sum_{R \subseteq Q_j}    \langle T(E^{b_1}_{Q_j} f),  \Delta^{b_2}_R g  \rangle_{\nu} \right|  \\
& \le  \sum_{j=1}^{2^n} \| 1_{Q_j} T(E^{b_1}_{Q_j} f) \|_{L^2(\nu)} \Big \|  \sum_{R \subseteq Q_j} \Delta^{b_2}_R g \Big \|_{L^2(\nu)} \\
&\lesssim \sum_{j =1}^{2^n} \left | \frac{\langle f \rangle^\mu_{Q_j}}{ \langle b_1 \rangle^\mu_{Q_j}}  \right |  \| 1_{Q_j} T(b_1 1_{Q_j} ) \|_{L^2(\nu)}  \| 1_{Q_j} g \|_{L^2(\nu)}  \\
& \lesssim \| g\|_{L^2(\nu)} \sum_{j =1}^{2^n} \mu(Q_j)^{-1/2} \left | \int_{Q_j} f d\mu \right| \\
& \lesssim \| g\|_{L^2(\nu)}  \| f \|_{L^2(\mu)},
\end{aligned}
\]
where we used \eqref{eqn:bd1}, testing condition (a), and H\"older's inequality. Clearly $S_3$ can be handled in an analogous manner.
 Similarly, if we consider $S_4$, testing condition (b) implies that
\[
\begin{aligned}
|S_4| &\le  \sum_{j,k =1}^{2^n}  | \langle T(E^{b_1}_{Q_j} f), E^{b_2}_{Q_k} g \rangle_{\nu} | \\
&\lesssim  \sum_{j,k =1}^{2^n} |\langle f\rangle_{Q_j}^{\mu} | |\langle g\rangle_{Q_k}^{\mu} | | \langle T(1_{Q_j} b_1), 1_{Q_k}b_2\rangle_{\nu}| \\
& \lesssim \sum_{j,k =1}^{2^n} |\langle f\rangle_{Q_j}^{\mu} | |\langle g\rangle_{Q_k}^{\mu} | \mu(Q_j)^{1/2} \nu(Q_k)^{1/2} \\
&\lesssim \| f\|_{L^2(\mu)} \| g\|_{L^2(\nu)}.
\end{aligned} \]
Now decompose $S_1$ as follows:
\[
\begin{aligned}
S_1 &= S_{11} + S_{12} + S_{13}+S_{14}\\
& =  \sum_{\substack{Q \subseteq \cup Q_j, R \subseteq \cup Q_k \\ \ell(R) <2^{-r} \ell(Q)}} \langle T(\Delta^{b_1}_Q f), \Delta^{b_2}_R g \rangle_{\nu}
+  \sum_{\substack{Q \subseteq \cup Q_j, R \subseteq \cup Q_k\\ \ell(Q) <2^{-r} \ell(R)}} \langle T(\Delta^{b_1}_Q f), \Delta^{b_2}_R g \rangle_{\nu}  \\
&+   \sum_{\substack{Q \subseteq \cup Q_j, R \subseteq \cup Q_k \\ 2^{-r} \ell(Q) \le \ell(R) \le \ell(Q)}}  \langle T(\Delta^{b_1}_Q f), \Delta^{b_2}_R g \rangle_{\nu} +
 \sum_{\substack{Q \subseteq \cup Q_j, R \subseteq \cup Q_k \\ 2^{-r} \ell(R) \le \ell(Q) < \ell(R)}}  \langle T(\Delta^{b_1}_Q f), \Delta^{b_2}_R g \rangle_{\nu}.
 \end{aligned}
\]
First we consider $S_{13}.$  Fix $Q\subseteq \cup Q_j$ and observe that if $R \subseteq \cup Q_k$ with $2^{-r} \ell(Q) \le \ell(R) \le \ell(Q)$, then the definition of well-localized implies that
\[   \langle T(\Delta^{b_1}_Q f), \Delta^{b_2}_R g \rangle_{\nu}\]
can only be nonzero if $R \subseteq Q^{(r)}.$ It is easy to show that there are only finitely many $R$ satisfying both $2^{-r} \ell(Q) \le \ell(R) \le  \ell(Q)$ and $R \subseteq Q^{(r)}$. Let $M_Q$ denote the number of such $R$, and label the $R$ cubes $R_Q^1, \dots, R_Q^{M_Q}$. Then $M_Q$ can be bounded by a constant $M$ that depends only on $n$ and $r$, not $Q$.
Similarly, one can show that
each $R$ can be an $R_{Q}^m$ for at most $N$ cubes $Q$, where $N$ is a constant depending on $n$ and $r$, but not on $R$.
Then using testing condition (b), Remark \ref{rem:MD}, and Lemma \ref{lem:decomp0}, we have
\[
\begin{aligned}
 |S_{13} | & \le \sum_{Q \subseteq \cup Q_j} \sum_{m=1}^{M_Q}  | \langle T(\Delta^{b_1}_Q f), \Delta^{b_2}_{R^m_Q} g \rangle_{\nu}  |  \\
 &  \lesssim \sum_{Q \subseteq \cup Q_j} \sum_{m=1}^{M_Q} \| \Delta^{b_1}_Q f \|_{L^2(\mu)}  \| \Delta^{b_2}_{R^m_Q} g \|_{L^2(\nu)}  \\
 & \le \left( \sum_{Q \subseteq \cup Q_j} \sum_{m=1}^{M_Q}  \| \Delta^{b_1}_Q f \|_{L^2(\mu)}^2 \right)^{1/2} \left( \sum_{Q \subseteq \cup Q_j} \sum_{m=1}^{M_Q} \| \Delta^{b_2}_{R^m_Q} g \|_{L^2(\nu)}^2 \right)^{1/2} \\
 & \lesssim  \| f\|_{L^2(\mu)}  \| g \|_{L^2(\nu)}.
\end{aligned}
\]
The sum $S_{14}$ can be handled in an analogous way.

Lastly, we consider $S_{11}$; by symmetry, the arguments given here, applied to $T^*$ instead of $T$, will also handle $S_{12}$. Observe that if $j \ne k$, then the definition of well-localized gives
\[ \sum_{\substack{Q \subseteq Q_j, R \subseteq Q_k \\ \ell(R) <2^{-r} \ell(Q)}} \langle T(\Delta^{b_1}_Q f), \Delta^{b_2}_R g \rangle_{\nu} =0\]
since for each for $Q' \in \ch Q$, $\ell(R) \le 2^{-r} \ell (Q')$ but $ R \not \subseteq Q'.$  Thus, we need only consider
\[ S_{11} =  \sum_{j=1}^{2^n} \sum_{\substack{Q,R \subseteq Q_j \\ \ell(R) <2^{-r} \ell(Q)}} \langle T(\Delta^{b_1}_Q f), \Delta^{b_2}_R g \rangle_{\nu}= \sum_{j=1}^{2^n} \sum_{\substack{R \subseteq Q \subseteq Q_j \\ \ell(R) <2^{-r} \ell(Q)}} \langle T(\Delta^{b_1}_Q f), \Delta^{b_2}_R g \rangle_{\nu},\]
where again we used the definition of well-localized.  This sum collapses as follows:
\[
\begin{aligned}
S_{11} &=  \sum_{j=1}^{2^n} \sum_{\substack{R \subseteq Q \subseteq Q_j \\ \ell(R) <2^{-r} \ell(Q)}} \left(  \sum_{Q' \in \ch Q} \langle T(E^{b_1}_{Q'} f), \Delta^{b_2}_R g \rangle_{\nu} - \langle T(E^{b_1}_{Q} f), \Delta^{b_2}_R g \rangle_{\nu} \right) \\
 &=   \sum_{j=1}^{2^n}  \sum_{\substack{ R \subseteq Q \subsetneq Q_j \\ \ell(R) = 2^{-r} \ell(Q)}}   \langle T(E^{b_1}_Q f), \Delta^{b_2}_R g \rangle_{\nu} -  \sum_{j=1}^{2^n} \sum_{\substack{ R \subseteq Q_j \\ \ell(R) < 2^{-r} \ell (Q_j)}}   \langle T(E^{b_1}_{Q_j} f), \Delta^{b_2}_R g \rangle_{\nu} \\
& : = S_{111} + S_{112}.
\end{aligned}
\]
To control the sum $S_{112}$, observe that using earlier arguments and \eqref{eqn:bd1}, we have
\[
\begin{aligned}
|S_{112}| &\le  \sum_{j=1}^{2^n} \left \| 1_{Q_j} T(E^{b_1}_{Q_j} f) \right \|_{L^2(\nu)} \Big \|  \sum_{\substack{ R \subseteq Q_j \\ \ell(R) < 2^{-r} \ell (Q_j)}}  \Delta^{b_2}_R g \Big \|_{L^2(\nu)} \\
& \lesssim \| f\|_{L^2(\mu)} \sum_{j=1}^{2^n} \left( \Big \| \sum_{\substack{ R \subseteq Q_j }}  \Delta^{b_2}_R g \Big\|_{L^2(\nu)}  +  \sum_{\substack{ R \subseteq Q_j \\ \ell(R) \ge 2^{-r} \ell (Q_j)}}  \|  \Delta^{b_2}_R g \|_{L^2(\nu)}\right) \\
& \lesssim \| f\|_{L^2(\mu)}  \| g\|_{L^{2}(\nu)},
\end{aligned}
\]
since there are at most finitely many terms in the last sum. Now we just need to consider $S_{111}$. Since $\Delta^{b_2}_R$ is a projection, we have
\[
\begin{aligned}
|S_{111}| &\lesssim \sum_{j=1}^{2^n}  \sum_{\substack{ R \subseteq Q \subsetneq Q_j \\ \ell(R) = 2^{-r} \ell(Q)}}  |\langle f \rangle^{\mu}_Q|  \left | \left  \langle (\Delta^{b_2}_R)^*  T(1_Q b_1 ), \Delta^{b_2}_R g \right \rangle_{\nu} \right |  \\
&\lesssim \sum_{j=1}^{2^n}  \sum_{\substack{ R \subseteq Q \subsetneq Q_j \\ \ell(R) = 2^{-r} \ell(Q)}}  |\langle f \rangle^{\mu}_Q| \| (\Delta^{b_2}_R)^*  T(1_Q b_1 ) \|_{L^2(\nu)} \| \Delta^{b_2}_R g  \|_{L^2(\nu)}  \\
&\le \sum_{j=1}^{2^n} \left(  \sum_{Q \subsetneq Q_j}   |\langle f \rangle^{\mu}_Q|^2  \sum_{\substack{ R \subseteq Q \\ \ell(R) = 2^{-r} \ell(Q)}} \| (\Delta^{b_2}_R)^*  T(1_Q b_1 ) \|^2_{L^2(\nu)} \right)^{1/2}  \left( \sum_{R \subseteq Q_j} \| \Delta^{b_2}_R g  \|^2_{L^2(\nu)} \right)^{1/2},
\end{aligned}
\]
where we used the Cauchy-Schwarz inequality. Lemma \ref{lem:decomp0} implies that the second term is bounded by $\|g\|_{L^2(\nu)}$. To control the first term, we need to show that the sequence $(a_Q)$, defined by
\[ a_Q := \sum_{\substack{ R \subseteq Q \\ \ell(R) = 2^{-r} \ell(Q)}} \| (\Delta^{b_2}_R)^*  T(1_Q b_1 ) \|^2_{L^2(\nu)}, \text{ for } Q \subsetneq Q_j,\]
and $a_Q:=0$ otherwise, is a $\mu$-Carleson sequence for each $j$.  Then the result follows by the dyadic Carleson embedding theorem, given in Theorem \ref{theo: Carleson embedding}.

To show that $(a_Q)$ is a $\mu$-Carleson sequence, we need only consider the case when $H \subsetneq Q_j$ for some $j$.
In particular, we need to control
\[ \sum_{Q \subseteq H} a_Q = \sum_{Q \subseteq H}  \sum_{\substack{ R \subseteq Q \\ \ell(R) = 2^{-r} \ell(Q)}} \| (\Delta^{b_2}_R)^*  T(1_Q b_1 ) \|^2_{L^2(\nu)}\]
by $\mu(H)$. To proceed, fix $Q \subseteq H$ and $R \subseteq Q$ with $\ell(R) = 2^{-r} \ell(Q)$. Then if $Q' \subseteq H$ with $Q \ne Q'$ and $\ell(Q) = \ell(Q')$, we immediately have $\ell(R) =2^{-r}\ell(Q')$ and $ R \not \subseteq Q'$. Then the definition of well-localized implies that
\[
\begin{aligned}
\| (\Delta^{b_2}_R)^*  T(1_{Q} b_1 ) \|^2_{L^2(\nu)} &=  \langle T(1_{Q}b_1), \Delta^{b_2}_R (\Delta^{b_2}_R)^*1_Q T(1_{Q}b_1 ) \rangle_{\nu} \\
&= \langle T(1_{H}b_1), \Delta^{b_2}_R (\Delta^{b_2}_R)^*T(1_{Q}b_1 ) \rangle_{\nu}  \\
&=  \langle  \Delta^{b_2}_R (\Delta^{b_2}_R)^* 1_H T(1_{H}b_1),T(1_{Q}b_1 ) \rangle_{\nu} \\
&=  \langle  \Delta^{b_2}_R (\Delta^{b_2}_R)^* 1_H T(1_{H}b_1),T(1_{H}b_1 ) \rangle_{\nu}  \\
& = \| (\Delta^{b_2}_R)^* 1_H T(1_{H} b_1 ) \|^2_{L^2(\nu)}. \end{aligned}
\]
Thus we can rewrite our sum as
\[ \sum_{Q \subseteq H} a_Q = \sum_{Q \subseteq H}  \sum_{\substack{ R \subseteq Q \\ \ell(R) = 2^{-r} \ell(Q)}} \| (\Delta^{b_2}_R)^* 1_H T(1_H b_1 ) \|^2_{L^2(\nu)}
 \lesssim \| 1_H T(1_H b_1) \|^2_{L^2(\nu)} \lesssim \mu(H),\]
where we used Lemma \ref{lem:decomp0} and testing condition (a). This shows that $(a_Q)$ is a $\mu$-Carleson sequence and completes the proof.
 \end{proofof}

\section{Local $Tb$ Theorem} \label{sec:LTB}

Before defining the system of test functions,  recall a standard notion of sparsity; a set $\mathcal{S} \subseteq \mathcal{D}$ is $\mu$-sparse if for all $R \in \mathcal{D}$,\[ \sum_{Q \in \mathcal{S}: Q\subseteq R} \mu(Q) \lesssim \mu(R).\]
Equivalently, the sequence $(a_Q)_{Q\in \mathcal{D}}$ defined by $a_Q = \mu(Q)$ for $Q \in \mathcal{S}$ and $a_Q = 0$ otherwise is a $\mu$-Carleson sequence.

\begin{defi} \label{def:sas} We say a system of functions $\{b_Q\}_{Q\in\mathcal{D}}$ is a \emph{sparse $L^{\infty}(\mu)$-accretive system} if it satisfies two conditions. First, $\{b_Q\}_{Q\in\mathcal{D}}$ is an \emph{$L^{\infty}(\mu)$-accretive system}, which here means that
\begin{itemize}
	\item[\rm (i)] $\supp (b_Q) \subseteq Q$
	\item[\rm (ii)] $\|b_Q\|_{L^\infty(\mu)} \lesssim 1$
	\item[\rm (iii)] $\big |\int_Q b_Q d \mu \big| \gtrsim \mu(Q)$,
\end{itemize}
for each $Q \in \mathcal{D}$, where the implied constants are independent of $Q$. Second, the set of cubes where the $b_Q$ change between generations is sparse. In particular, if $\mathcal{S}_b := \{ Q\in \mathcal{D}: b_Q \ne b_{Q^{(1)}} 1_Q \}$, then $\mathcal{S}_b$ is $\mu$-sparse.
\end{defi}

The definition of an $L^{\infty}(\mu)$-accretive system given above is very similar to the definitions used in both \cite{HM2012, NTV2002}, but does not impose conditions on any $\| T b_Q\|_{L^{\infty}(\mu)}$. The testing conditions we use appear later. Then given a sparse $L^{\infty}(\mu)$-accretive system $\{b_Q\}_{Q\in\mathcal{D}}$, we can partition $\mathcal{D}$ into two sets: $\mathcal{D}_b$ and $\mathcal{C}_b$. $\mathcal{D}_b$ will denote the set of $Q \in \mathcal{D}$ that are contained in some $P \in \mathcal{S}_b$. The minimal such $P$ will be denoted by $P_Q$. Similarly $\mathcal{C}_b = \mathcal{D} \setminus \mathcal{D}_b$ will denote the set of cubes that are not contained in any $P \in \mathcal{S}_b.$ Note that if a point $x$ is in two cubes $Q,R \in \mathcal{C}_b,$ then $b_Q(x) = b_R(x)$. This means that if we set
\begin{equation} \label{eqn:bcb} b(x) := b_Q(x) \text{ for each  $x$ in some $Q \in \mathcal{C}_b$}\end{equation}
and $b(x):=0$ otherwise, then $b$ is well defined on $\mathbb{R}^n$ and satisfies $b_Q = b 1_Q$ for every $Q \in \mathcal{C}_b.$

\begin{rema} If $\{\tilde{b}_Q\}_{Q\in \mathcal{D}}$ is an $L^{\infty}(\mu)$-accretive system and if $\mu$ is compactly supported, then it can be used to create a sparse $L^{\infty}(\mu)$-accretive system. This is basically the stopping-time construction from \cite[pp. 4823]{HM2012} and \cite[pp. 269]{NTV2002}. In what follows, without loss of generality, we assume the implied constant in property (iii) of Definition \ref{def:sas} is some positive $\delta <1$ and the implied constant in (ii) of Definition \ref{def:sas} is some $C>1.$

First choose cubes $\{ Q_0^j\}_{j=1}^{2^n}$ with no common ancestors such that $\supp \mu \subseteq \cup_j Q_0^j.$  Then set $\mathcal{D}^0 = \{Q^j_0\}_j$ and for each $Q^j_0$, collect all maximal cubes $Q \subsetneq Q^j_0$ satisfying
	$| \int_Q \tilde{b}_{Q^j_0} d \mu|  < \delta^2 \mu(Q).
	$
Denote the resulting collection of cubes by $\mathcal{D}^1$. Then for each cube $Q_1^k \in \mathcal{D}^1$, collect all the maximal cubes $Q \subsetneq Q^k_{1}$ satisfying
	$| \int_Q \tilde{b}_{Q^k_1} d \mu|  < \delta^2 \mu(Q),
	$
and denote the resulting collection by $\mathcal{D}^2.$ Proceeding in this manner gives collections $\mathcal{D}^j$ for every $j \in \mathbb{N}$.
Using arguments appearing in \cite{NTV2002}, for $\tau = \frac{C-\delta}{C-\delta^2}<1$ and $R \in \mathcal{D}^j$, one can show
\[\sum_{Q \in \mathcal{D}^{j+1}: Q \subseteq R} \mu(Q) =  \mu \left( \cup_{Q \in \mathcal{D}^{j+1}} Q \cap R \right)  \lesssim \tau \mu(R),\]
where the implied constant does not depend on $R$. Then a simple argument shows that  these stopping cubes are $\mu$-sparse, namely for all $R\in \mathcal{D}$,
  \begin{equation} \label{eqn:sparse}
    \sum_{Q \in \cup_j \mathcal{D}^j, Q \subseteq R} \mu(Q) \lesssim \mu(R).
    \end{equation}

We can define the associated sparse $L^{\infty}(\mu)$-accretive system as follows. First for $Q$ with $Q \cap( \cup_j Q_0^j) = \emptyset$, let $b_Q \equiv 0.$  For $Q$ with $Q^j_0 \subseteq Q$ for some $j$, let $b_Q = \tilde{b}_{Q^j_0}.$ Then these $b_Q$ trivially satisfy (i)-(iii) in Definition \ref{def:sas}.
For each $Q \subseteq \cup_j Q_0^j$, let $Q^a$ denote the smallest cube in $\cup_j \mathcal{D}^j$ containing $Q$ and set $b_Q := \tilde{b}_{Q_a} 1_Q$. It is easy to check that these $b_Q$ also satisfy conditions (i)-(iii). Conditions (i) and (ii) are immediate. Similarly, if  $Q^a=Q$, condition (iii) follows.   If $Q \ne Q^a$, then by construction, $ | \int_Q \tilde{b}_{Q^a} d \mu| \ge \delta^{2} \mu(Q).$ Moreover, $b_Q \ne b_{Q^{(1)}}\text{1}_Q$ implies that $Q \in \cup_j \mathcal{D}^j$. Thus if we define $\mathcal{S}_b$ as in Definition \ref{def:sas}, then \eqref{eqn:sparse} implies that $\mathcal{S}_b \subseteq  \cup_j \mathcal{D}^j$ is $\mu$-sparse, as needed.
\end{rema}

Let $\{b_Q\}$ be a sparse $L^{\infty}(\mu)$-accretive system. Then the functions in $L^2(\mu)$ can be decomposed using these accretive systems. First define the associated expectations and martingale differences
	\[
	E^\mu_Q f = \frac{\langle f\rangle^\mu_Q}{\langle b_{Q}\rangle^\mu_Q} b_{Q}
	\ \text{ and } \
	\Delta^\mu_Q f = \sum_{Q' \in \ch Q} E^\mu_{Q'} f - E^\mu_Q f
	\]
	for all $Q \in \mathcal{D}.$ It is worth pointing out that, to make the two setups easier to differentiate, this notation $E^\mu_Q$ and $\Delta^\mu_Q$ is different from the notation $E^b_Q$ and $\Delta^b_Q$ in Sections~\ref{sec:GTB} and~\ref{sec:GTBP}. Now note that each $\Delta^\mu_Q f$ is supported on $Q$ and satisfies $\langle  \Delta^\mu_Q f \rangle_Q^{\mu} = 0.$ Because of this, we call all functions in the range space  $\Delta^\mu_Q L^2(\mu)$ $\mu$-Haar functions associated to $Q$ and will denote these functions by $h^\mu_Q$. Unlike the classical situation, the spaces $\Delta^\mu_Q L^2(\mu)$ and $\Delta^\mu_R L^2(\mu)$ need not be orthogonal for $Q \ne R.$ One  can also compute
\begin{equation} \label{eqn:dual} (\Delta^\mu_Q)^* f  = \sum_{Q' \in \ch Q} \left( \frac{ \langle b_{Q'} f \rangle^{\mu}_{Q'}}{\langle b_{Q'} \rangle^\mu_{Q'}} - \frac{\langle b_Q f\rangle_Q^{\mu}}{ \langle b_Q\rangle_Q^{\mu}} \right) 1_{Q'},\end{equation}
for each $Q\in \mathcal{D}$.
The arguments in \cite[pp. 4824-4825]{HM2012} and \cite[pp. 271-274]{NTV2002} adapt to this setting to give the decomposition below and testing condition (i). Because our setup is somewhat different and the details for (ii) do not appear in \cite{HM2012, NTV2002}, we give the proof of the following lemma in the appendix.

\begin{lem} \label{lem:decomp}
Let $\{b_Q\}$ be a sparse $L^{\infty}(\mu)$-accretive system and let $f \in L^2(\mu)$. Then for each $d\in \mathbb{Z}$,
\[ f = \sum_{\substack{ Q \in \mathcal{D} \\ \ell(Q) \le 2^d}} \Delta^\mu_Q f + \sum_{\substack{Q \in \mathcal{D} \\ \ell(Q) = 2^d}} E^\mu_Q f\]
with convergence in $L^2(\mu).$ Moreover, the following estimates hold
\begin{itemize}
\item[(i)] $\displaystyle \sum_{\substack{ Q \in \mathcal{D} \\ \ell(Q) \le 2^d}} \| \Delta^\mu_Q f \|^2_{L^2(\mu)}   + \sum_{\substack{Q \in \mathcal{D} \\ \ell(Q) = 2^d}}  \| E^\mu_{Q} f \|^2_{L^2(\mu)}  \lesssim \|f\|_{L^2(\mu)}^2,$
\item [(ii)] $\displaystyle \sum_{\substack{Q \in \mathcal{D}}}  \| (\Delta^\mu_Q)^* f \|^2_{L^2(\mu)} \lesssim \|f\|_{L^2(\mu)}^2.$
\end{itemize}
\end{lem}

A simple estimate gives that each $\| E^\mu_{Q} f \|_{L^2(\mu)} \lesssim \| f 1_{Q} \|_{L^2(\mu)}$, and thus
\begin{equation} \label{eqn:bd} \Big \|  \sum_{\substack{ Q \in \mathcal{D} \\ \ell(Q) \le 2^d}}  \Delta^\mu_Q f \Big \|_{L^2(\mu)} \lesssim \|f\|_{L^2(\mu)}.\end{equation}
To see how much these $\Delta_Q^\mu$ differ from projections, one can compute
	\begin{equation} \label{eqn:error}
	\Delta^\mu_Q f  - (\Delta^\mu_Q)^2 f =  \sum_{\substack{P \in \ch Q \cap \mathcal{S}_{b}}} \varphi^\mu_P, \ \text{ where } \
	\varphi^\mu_P = \frac{\langle f\rangle^\mu_Q}{\langle b_{Q}\rangle^\mu_Q} \left(\frac{\langle b_{Q}\rangle^\mu_P}{\langle b_P\rangle^\mu_P} b_P - b_{Q}\right) 1_P.
	\end{equation}
Then properties (i)-(iii) of $\{b_Q\}_{Q\in \mathcal{D}}$ imply that $\|\varphi_P^{\mu}\|_{L^2(\mu)} \lesssim  | \langle f\rangle^\mu_Q | \mu(P)^{1/2}.$

In what follows, we will examine pairs of sparse accretive systems associated to two Borel measures.

\begin{defi} \label{def:sas2} We say a system of functions $\textbf{b} = \{b^1_Q, b^2_Q\}_{Q\in\mathcal{D}}$ is a \emph{sparse $L^{\infty}(\mu,\nu)$-accretive system} if $\{b^1_Q\}$ is a sparse $L^{\infty}(\mu)$-accretive system, $\{b^2_Q\}$ is a sparse $L^{\infty}(\nu)$-accretive system, and this additional sparsity condition holds: $\mathcal{S}_{b_1}$ is sparse with respect to $\nu$ and $\mathcal{S}_{b_2}$ is sparse with respect to $\mu.$
\end{defi}

\begin{rema} The extra sparsity condition in Definition \ref{def:sas2} implies that the set of cubes where the $b^j_Q$ change between generations is small with respect to both measures. Trivially, this condition is satisfied if as in Section \ref{sec:GTB}, for $j=1,2$, there is one $b^j$ so that $b^j_Q = b^j 1_Q$ for all $Q \in \mathcal{D}$. Similarly this condition is satisfied if $\mu =\nu.$ So this setup generalizes both the accretive function case and the one-weight case.
\end{rema}

Now let $\textbf{b} = \{b^1_Q, b^2_Q\}_{Q\in\mathcal{D}}$ be a sparse $L^{\infty}(\mu,\nu)$-accretive system. We say $T$ is an operator acting formally from $L^2(\mu)$ to $L^2(\nu)$ with respect to $\textbf{b}$ if its bilinear form $\left\langle Tb^1_Q, b^2_R \right\rangle_\nu$  is well defined for all $Q,R \in \mathcal{D}$.  Then we can define the well-localized operators in this setting.

\begin{defi}\label{defi: WLO}
Let $\textbf{b}$ be a sparse $L^{\infty}(\mu,\nu)$-accretive system and let $T$ be an operator acting formally from $L^2(\mu)$ to $L^2(\nu)$ with respect to $\textbf{b}.$
We say that $T$ is \emph{lower triangularly localized with respect to $\textbf{b}$ with radius $r$} if there exists an integer $r \ge 0$ such that for all cubes $Q,R \in \mathcal{D}$ with $\ell(R) \leq 2 \ell(Q)$ and all $\nu$-Haar functions $h^\nu_R$ on $R$
	\[
	\left\langle T (b^1_{Q}), h^\nu_R \right\rangle_\nu = 0
	\]
if $R \not \subseteq Q^{(r)}$ or if $\ell(R) \leq 2^{-r}\ell(Q)$ and $R \not\subseteq Q$.
We say that the operator $T$ is \emph{well-localized with respect to~\textbf{b} of radius $r$} if both $T$ and its formal adjoint $T^*$ are lower triangularly localized with respect to $\textbf{b}$ with radius $r$ and if $T$ (and $T^*$) satisfy an additional localization property: for $T$, if $Q \subseteq S$ with $P^1_Q =P^1_S$ (here $P^1_Q$ is the minimal $P \in \mathcal{S}_{b_1}$ with $Q \subseteq P$) or if both $S, Q \in \mathcal{C}_{b_1}$, then for $R \in \ch^{r}(Q)$
\begin{equation} \label{eqn:csc} \| (\Delta^{\nu}_R)^* T(b^1_Q) \|_{L^2(\nu)}^2= \| (\Delta^{\nu}_R)^* T(b^1_S 1_Q) \|_{L^2(\nu)}^2  \lesssim \| (\Delta^{\nu}_R)^* T(b^1_S) \|_{L^2(\nu)}^2.\end{equation}
\end{defi}

\begin{rema} Condition \eqref{eqn:csc} is a new and somewhat restrictive condition that we need for the proof to work. If possible, we would like to relax this condition so that the theorem applies to more operators. However, in the accretive function setting with $\{b_1, b_2\}$ as in Section~\ref{sec:GTB}, this condition follows immediately from the other parts of the well-localized definition. Indeed, in the context of Theorem \ref{theo: Tb-Theorem}, Condition \eqref{eqn:csc} is  trivial whenever each $Q' \subseteq S$ with $\ell(Q') = \ell(Q)$ also satisfies $P^1_{Q'} = P^1_S$ or each satisfies $Q' \in \mathcal{C}_{b_1}$ respectively. To see this, note that in those cases, for $R \in \ch^r(Q)$,
\[  \| (\Delta^{\nu}_R)^* T(b^1_S 1_{Q'}) \|_{L^2(\nu)}^2  =  \| (\Delta^{\nu}_R)^* T(b^1_{Q'}) \|_{L^2(\nu)}^2 =\left \langle  T(b^1_{Q'}),  \Delta^\nu_R (\Delta^{\nu}_R)^*  T(b^1_{Q'})\right \rangle_{\nu} =0\]
because $\ell(R) = 2^{-r}\ell(Q')$ and $R \not\subseteq Q'$. Then it is immediate that
\[  \| (\Delta^{\nu}_R)^* T(b^1_Q) \|_{L^2(\nu)}^2= \Big \|  \sum_{Q' \subseteq S, \ell(Q) = \ell(Q')} (\Delta^{\nu}_R)^* T(b^1_S 1_{Q'}) \Big \|_{L^2(\nu)}^2 =  \| (\Delta^{\nu}_R)^* T(b^1_S) \|_{L^2(\nu)}^2,\]
as needed.
\end{rema}
Then we can prove the following local $Tb$ theorem.

\begin{theo}\label{theo: Tb-Theorem} Let $T$ be a well-localized operator with respect to a sparse $L^{\infty}(\mu,\nu)$-accretive system \textbf{b} with radius $r$. Further assume
\begin{itemize}
\item[(a)] $	\|T(b^1_{Q})\|^2_{L^2(\nu)} \lesssim \mu(Q)  \ \text{ and } \
	\|T^*(b^2_{Q})\|^2_{L^2(\mu)} \lesssim \nu(Q), \quad \forall Q\in\mathcal{D}$;
\item[(b)] For all $Q,R \in \mathcal{D}$ satisfying $2^{-r}\ell(Q) \leq \ell(R) \leq 2^r \ell(Q)$,
\begin{equation}\label{eq: T(b)-Theo_cond2}
	\left|\left\langle T\Delta^\mu_Q f, \Delta^\nu_R g\right\rangle_\nu\right|
	\lesssim \|\Delta^\mu_Q f\|_{L^2(\mu)}\|\Delta^\nu_R g\|_{L^2(\nu)}
	\quad \forall f \in L^2(\mu), g \in L^2(\nu);
\end{equation}
\item[(c)] For all $Q \in \mathcal{D}$ and $P \in \ch^{r+1}(Q)$,
\[
P \in \mathcal{S}_{b_2} \text{ implies } \left \| 1_P T (b^1_{Q})\right \|^2_{L^2(\nu)} \lesssim \mu(P)
\text{ and } P \in \mathcal{S}_{b_1} \text{ implies }  \left \| 1_P T^* (b^2_{Q})\right \|^2_{L^2(\mu)} \lesssim \nu(P).
 \]
\end{itemize}
Then $T: L^2(\mu) \to L^2(\nu)$ is bounded.
\end{theo}

\begin{rema} A couple remarks about the testing conditions are in order. First, conditions (a) and (b) are similar to, but somewhat different than, the testing conditions in Theorem~\ref{thm:GTB}. However, if our operator $T$ is further localized in the sense that
\[ \langle T (b^1_Q), b^2_R \rangle_{\nu} =0,\]
if $Q, R \in \mathcal{D}$ have no common ancestors, then we can replace this testing condition (a) with the condition from Theorem \ref{thm:GTB}:
\[ \| 1_QT(b^1_{Q})\|^2_{L^2(\nu)} \lesssim \mu(Q)  \ \text{ and } \
	\|1_Q T^*(b^2_{Q})\|^2_{L^2(\mu)} \lesssim \nu(Q), \quad \forall Q\in\mathcal{D}.\]
Condition (b) is necessarily different in this setting because the martingale differences are more complicated for accretive systems.

Meanwhile, testing condition (c) did not appear in Theorem \ref{thm:GTB}. Indeed, in the case of accretive functions $\{b_1, b_2\}$, condition (c) is trivial because $\mathcal{S}_{b_1}, \mathcal{S}_{b_2} = \emptyset.$ Similarly, if $\nu$ and $\mu$ are doubling measures, then (c) is immediate. To see this, note that because $P \in \ch^{r+1}(Q)$, the doubling condition implies that $P$ and $Q$ have comparable $\mu$-sizes. Then  testing condition (a) immediately implies
\[  \left \| 1_P T (b^1_{Q})\right \|^2_{L^2(\nu)}  \lesssim \mu(Q) \lesssim \mu(P),\]
and a similar argument controls $ \| 1_P T^* (b^2_{Q}) \|^2_{L^2(\mu)}.$
\end{rema}

\section{Proof of Theorem \ref{theo: Tb-Theorem}} \label{sec:LTBP}


Now let us consider the proof of Theorem \ref{theo: Tb-Theorem}:\\

\noindent \begin{proofof}{Theorem~\ref{theo: Tb-Theorem}} Fix $f \in L^2(\mu)$ and $g \in L^2(\nu)$ and without loss of generality, assume that they are compactly supported. Then, there is an integer $d$ and cubes $Q_1, \dots, Q_{2^n} \in \mathcal{D}$ with no common ancestors such that $\ell(Q_j) = 2^d$ and $\supp f, \supp g \subseteq \cup Q_j.$ By Lemma \ref{lem:decomp}, we can write
\[
\begin{aligned}
f & = \sum_{Q \subseteq \cup Q_j} \Delta^{\mu}_Q f + \sum_{j=1}^{2^n} E^{\mu}_{Q_j} f = f_1 + f_2 \\
g &= \sum_{R \subseteq \cup Q_k} \Delta^{\nu}_R g + \sum_{k=1}^{2^n} E^{\nu}_{Q_k}  g = g_1 + g_2.
\end{aligned}
\]
By duality, it suffices to show that $|\langle Tf,g\rangle_{\nu}| \le C \| f\|_{L^2(\mu)} \| g \|_{L^2(\nu)}.$ We break the inner product into the following four terms
\[ S_1 + S_2 +S_3 +S_4 = \langle Tf_1,g_1\rangle_{\nu} + \langle Tf_2,g_1\rangle_{\nu} + \langle Tf_1,g_2\rangle_{\nu} + \langle Tf_2,g_2\rangle_{\nu}\]
to handle separately. The sums $S_2$, $S_3$, and $S_4$ are handled in a way analogous to those in the proof of Theorem \ref{thm:GTB}, so we leave the details to the reader.


Now decompose $S_1$ as
\[
\begin{aligned}
S_1 &= S_{11} + S_{12}+S_{13}+S_{14} \\
& =  \sum_{\substack{Q \subseteq \cup Q_j, R \subseteq \cup Q_k \\ \ell(R) <2^{-r} \ell(Q)}} \langle T(\Delta^{\mu}_Q f), \Delta^{\nu}_R g \rangle_{\nu}
+  \sum_{\substack{Q \subseteq \cup Q_j, R \subseteq \cup Q_k\\ \ell(Q) <2^{-r} \ell(R)}} \langle T(\Delta^{\mu}_Q f), \Delta^{\nu}_R g \rangle_{\nu}  \\
&+   \sum_{\substack{Q \subseteq \cup Q_j, R \subseteq \cup Q_k \\ 2^{-r} \ell(Q) \le \ell(R) \le \ell(Q)}}  \langle T(\Delta^{\mu}_Q f), \Delta^{\nu}_R g \rangle_{\nu}+ \sum_{\substack{Q \subseteq \cup Q_j, R \subseteq \cup Q_k \\ 2^{-r} \ell(R) \le \ell(Q) <\ell(R)}}  \langle T(\Delta^{\mu}_Q f), \Delta^{\nu}_R g \rangle_{\nu}.
\end{aligned}
\]
As $S_{13}$ (and $S_{14}$) can be controlled as in the proof of Theorem \ref{thm:GTB}, we omit the details.  The main differences are using testing condition (b) and Lemma \ref{lem:decomp}.

Lastly consider sums $S_{11}$ and $S_{12}$. By symmetry, we need only estimate $S_{11}$. First, the definition of well-localized implies that when $ j \ne k$ the interior sums in $S_{11}$ vanish. Thus, we have
\begin{align*} S_{11} & = \sum_{j=1}^{2^n}
	\sum_{\substack{Q, R \subseteq Q_j \\ \ell(R) < 2^{-r}\ell(Q)}}  \left( \sum_{Q' \in \ch Q} \left\langle T(E^\mu_{Q'} f ), \Delta^\nu_{R} g\right\rangle_\nu  - \left\langle T(E^\mu_{Q} f), \Delta^\nu_{R} g\right\rangle_\nu  \right) \\
		&= \sum_{j=1}^{2^n} \sum_{\substack{R \subseteq Q \subsetneq Q_j \\ \ell(R) \leq 2^{-r}\ell(Q)}} \left\langle T(E^\mu_{Q} f ), \Delta^\nu_{R} g\right\rangle_\nu - \sum_{j=1}^{2^n} \sum_{\substack{R \subseteq Q \subseteq Q_j \\ \ell(R) < 2^{-r}\ell(Q)}}  \left\langle T(E^\mu_{Q} f), \Delta^\nu_{R} g\right\rangle_\nu \\
		&= \sum_{j=1}^{2^n}  \sum_{\substack{R \subseteq Q \subsetneq Q_j \\ \ell(R) = 2^{-r}\ell(Q)}}  \left\langle T(E^\mu_{Q} f), \Delta^\nu_{R} g\right\rangle_\nu- \sum_{j=1}^{2^n}  \sum_{\substack{R \subseteq Q_j \\ \ell(R) < 2^{-r}\ell(Q_j)}}  \left\langle T(E^\mu_{Q_j} f), \Delta^\nu_{R} g\right\rangle_\nu \\
		&:= S_{111} +S_{112},
\end{align*}
where the second equality used the definition lower triangularly localized. We can estimate $|S_{112}|$ easily by
	\[
	\begin{aligned}
	|S_{112}| &\lesssim \sum_{j=1}^{2^n} |\langle f\rangle^\mu_{Q_j}|  \|1_{Q_j} T(b^1_{Q_j})\|_{L^2(\nu)}  \left ( \Big \| \sum_{\substack{R \subset Q_j}}   \Delta^\nu_{R} g \Big \|_{L^2(\nu)} + \sum_{\substack{R\subseteq Q_j \\ \ell(R) \ge 2^{-r} \ell(Q_j)}} \|  \Delta^\nu_{R} g  \|_{L^2(\nu)} \right) \\
	&\lesssim \|f\|_{L^2(\mu)}\|g\|_{L^2(\nu)},
	\end{aligned}
	\]
	where the first $\Delta^\nu_R g$ sum is controlled using \eqref{eqn:bd} and the second sum is bounded because it only includes a finite number of terms.
One can now control $S_{111}$ by fixing $j$ and controlling
	\[
	\sum_{\substack{R \subseteq Q  \subsetneq Q_j \\ \ell(R) = 2^{-r}\ell(Q)}} \left|\frac{\langle f\rangle^\mu_Q}{\langle b^1_{Q}\rangle^\mu_Q} \left\langle T(b^1_{Q}), \Delta^\nu_{R} g\right\rangle_\nu\right|
	\lesssim T_1 + T_2,
	\]
where
	\[
	T_1 =\sum_{\substack{R \subseteq Q \subsetneq Q_j \\ \ell(R) = 2^{-r}\ell(Q)}} | \langle f\rangle^\mu_Q  \left\langle (\Delta^\nu_R)^*T (b^1_{Q}), \Delta^\nu_{R} g\right\rangle_\nu|
	\]
and
	\[
	T_2 = \sum_{\substack{R \subseteq Q \subsetneq Q_j \\  \ell(R) = 2^{-r}\ell(Q)}} \sum_{\substack{P \in \ch R \\ P \in \mathcal{S}_{b_2}}} | \langle f\rangle^\mu_Q  \left\langle T (b^1_{Q}), \varphi^\nu_P \right\rangle_\nu |,
	\]
	where $\varphi^\nu_P$ depends on $g$ and is defined in \eqref{eqn:error}.
By the Cauchy-Schwarz inequality and Lemma~\ref{lem:decomp}, we have
\begin{align*}
	T_1 &\leq  \Bigg( \sum_{\substack{R \subseteq Q \subsetneq Q_j \\ \ell(R) = 2^{-r}\ell(Q)}} \left| \langle f\rangle^\mu_Q \right|^2 \left\|(\Delta^\nu_R)^*T(b^1_{Q})\right\|^2_{L^2(\nu)}\Bigg)^{1/2}
	\Bigg(\sum_{R \subseteq Q_j} \left\|\Delta^\nu_{R} g\right\|^2_{L^2(\nu)}\Bigg)^{1/2} \\
		&\lesssim \Bigg( \sum_{Q \subsetneq Q_j}\left|\langle f\rangle^\mu_Q\right|^2 a_Q\Bigg)^{1/2} \|g\|_{L^2(\nu)},
\end{align*}
where
	$$
	a_Q = \sum_{R \in \ch^r Q} \left\|(\Delta^\nu_R)^*T(b^1_{Q})\right\|^2_{L^2(\nu)}
	$$
	for $Q \subsetneq Q_j$ and $a_Q=0$ otherwise. Then to apply the Carleson embedding theorem, we need to show $(a_Q)$ is a $\mu$-Carleson sequence. To do this, fix
	 $H\in\mathcal{D}$ and without loss of generality, assume $H \subsetneq Q_j$.  For now, assume $H \in \mathcal{D}_{b_1}$. This means there is some $P \in \mathcal{S}_{b_1}$ with $H \subseteq P.$ The minimal such $P$ is denoted by $P^1_H$ and for a general $Q \in  \mathcal{D}_{b_1}$, it is denoted $P^1_Q$.
	 Then we can write
	\[
	\begin{aligned}
	\sum_{Q \subseteq H} a_Q &= \sum_{\substack{Q \subseteq H  \\ P^1_Q = P^1_H}} \sum_{R \in \ch^r Q}  \left\|(\Delta^\nu_R)^*T(b^1_{Q})\right\|^2_{L^2(\nu)}
	+  \sum_{\substack{P \subsetneq H \\ P \in \mathcal{S}_{b_1}}} \sum_{\substack{Q \subseteq H \\ P^1_Q = P}} \sum_{R \in \ch^r Q}  \left\|(\Delta^\nu_R)^*T(b^1_{Q})\right\|^2_{L^2(\nu)}\\
 &:= T_{11} +T_{12}.
	\end{aligned}
	\]
We can control $T_{11}$ using the localization condition \eqref{eqn:csc} in the definition of well-localized,  the dual square function estimate in Lemma \ref{lem:decomp},  and testing condition (a) as follows:
\[
 T_{11}   \lesssim \sum_{\substack{Q \subseteq H  \\ P^1_Q = P^1_H}} \sum_{R \in \ch^r Q}  \left\|(\Delta^\nu_R)^*T(b^1_{H})\right\|^2_{L^2(\nu)}
   \lesssim   \left\| T(b^1_{H})\right\|^2_{L^2(\nu)}   \lesssim \mu(H),
\]
as needed.  The same arguments allow us to control $T_{12}$ as follows:
\[
\begin{aligned}
\sum_{\substack{P \subsetneq H \\ P \in \mathcal{S}_{b_1}}} \sum_{\substack{Q \subseteq H  \\ P^1_Q = P}} \sum_{R \in \ch^r Q}  \left\|(\Delta^\nu_R)^*T(b^1_{Q})\right\|^2_{L^2(\nu)} &\lesssim
\sum_{\substack{P \subsetneq H \\ P \in \mathcal{S}_{b_1}}} \sum_{\substack{Q \subseteq H  \\ P^1_Q = P}} \sum_{R \in \ch^r Q}  \left\|(\Delta^\nu_R)^*T(b^1_{P})\right\|^2_{L^2(\nu)} \\
  & \lesssim  \sum_{\substack{P \subsetneq H \\ P \in \mathcal{S}_{b_1}}}  \left\| T(b^1_{P})\right\|^2_{L^2(\nu)}  \\
  &  \lesssim  \sum_{\substack{P \subsetneq H \\ P \in \mathcal{S}_{b_1}}}  \mu(P) \lesssim \mu(H),
\end{aligned}
\]
where we used the fact that $\mathcal{S}_{b_1}$ is $\mu$-sparse. Now if $H$ was in $\mathcal{C}_{b_1}$, instead of $\mathcal{D}_{b_1}$, then $T_{12}$ would be the same, and $T_{11}$ would become
\[  \sum_{\substack{Q \subseteq H:  \\ Q \in \mathcal{C}_{b_1}}} \sum_{R \in \ch^r Q}  \left\|(\Delta^\nu_R)^*T(b^1_{Q})\right\|^2_{L^2(\nu)} \lesssim
  \sum_{\substack{Q \subseteq H:  \\ Q \in \mathcal{C}_{b_1}}} \sum_{R \in \ch^r Q}  \left\|(\Delta^\nu_R)^*T(b^1_{H})\right\|^2_{L^2(\nu)} \lesssim  \left\|T(b^1_{H})\right\|^2_{L^2(\nu)},  \]
so the same bound holds. Thus $(a_Q)$ is a $\mu$-Carleson sequence, so an application of the Carleson embedding theorem gives the bound for $T_1$. To control $T_2$, begin as follows:
\begin{align*}
	&T_2 = \sum_{\substack{R,Q: Q \subsetneq Q_j \\  R \in \ch^r(Q) }} \sum_{\substack{P \in \ch R \\ P \in \mathcal{S}_{b_2}}} | \langle f\rangle^\mu_Q| \left|\left\langle T (b^1_{Q}), \varphi^\nu_P \right\rangle_\nu \right|  \\
	& \lesssim \Big(  \sum_{\substack{R,Q: Q \subsetneq Q_j \\  R \in \ch^r(Q) }} \sum_{\substack{P \in \ch R \\ P \in \mathcal{S}_{b_2}}} | \langle f\rangle^\mu_Q|^2 \left \| 1_P T (b^1_{Q})\right \|^2_{L^2(\nu)} \Big)^{1/2} \Big(  \sum_{\substack{R,Q: Q \subsetneq Q_j \\  R \in \ch^r(Q) }} \sum_{\substack{P \in \ch R \\ P \in \mathcal{S}_{b_2}}} \left\| \varphi^\nu_P \right \|^2_{L^2(\nu)} \Big)^{1/2} \\
	& \lesssim \Big(  \sum_{Q \subsetneq Q_j} | \langle f\rangle^\mu_Q|^2 \sum_{\substack{P \in \ch^{(r+1)} Q \\ P \in \mathcal{S}_{b_2}}}  \left \| 1_P T (b^1_{Q})\right \|^2_{L^2(\nu)} \Big)^{1/2} \Big(  \sum_{R  \subsetneq Q_j} \sum_{\substack{P \in \ch R \\ P \in \mathcal{S}_{b_2}}} \left\| \varphi^\nu_P \right \|^2_{L^2(\nu)} \Big)^{1/2} \\
	&\lesssim \Big(  \sum_{Q \subsetneq Q_j} | \langle f\rangle^\mu_Q|^2 \sum_{\substack{P \in \ch^{(r+1)} Q \\ P \in \mathcal{S}_{b_2}}} \left \| 1_P T (b^1_{Q})\right \|^2_{L^2(\nu)} \Big)^{1/2}
	\Big(\sum_{\substack{R \subsetneq Q_j}} \left|\langle g\rangle^\nu_R\right|^2 \sum_{\substack{P \in \ch R \\ P \in \mathcal{S}_{b_2}}} \nu(P)\Big)^{1/2} \\
	&\lesssim \Big(\sum_{Q  \subsetneq Q_j} \left|\langle f\rangle^\mu_Q\right|^2 b_Q\Big)^{1/2} \|g\|_{L^2(\nu)},
\end{align*}
where
	$$
	b_Q = \sum_{\substack{P \in \ch^{(r+1)} Q \\ P \in \mathcal{S}_{b_2}}} \left \| 1_P T (b^1_{Q})\right \|^2_{L^2(\nu)},
	$$
for $Q \subsetneq Q_j$ and $b_Q =0$ otherwise. In the above computation, we also used the Carleson embedding theorem and the fact that $\mathcal{S}_{b_2}$ is $\nu$-sparse.
To complete the proof, we need to show that $(b_Q)$ is a $\mu$-Carleson sequence.
To do this, fix $H \in \mathcal{D}$ and without loss of generality, assume $H \subsetneq Q_j$. Then by testing condition (c), we have
\[ \sum_{Q \subseteq H} b_Q = \sum_{Q \subseteq H} \sum_{\substack{P \in \ch^{(r+1)} Q \\ P \in \mathcal{S}_{b_2}}} \left \| 1_P T (b^1_{Q})\right \|^2_{L^2(\nu)}
	\lesssim \sum_{\substack{ P \subseteq H \\ P \in \mathcal{S}_{b_2}}} \mu(P) \lesssim \mu(H), \]
	where we used the fact that $\mathcal{S}_{b_2}$ is also $\mu$-sparse.
\end{proofof}

\section{Appendix: Proof of Lemma \ref{lem:decomp}}
The proof requires the following well-known square function bound:

\begin{theo} \label{thm:USF} If $\mu$ is a Borel measure on $\mathbb{R}^n$ and $f \in L^2(\mu)$, then
\[ \sum_{Q \in \mathcal{D}} \left | \langle f \rangle_{Q^{(1)}}^\mu -  \langle f \rangle_{Q}^\mu \right |^2 \mu(Q) \lesssim \| f \|^2_{L^2(\mu)}.\]
\end{theo}

Let us proceed to the proof of Lemma \ref{lem:decomp}. \\

\noindent \begin{proofof}{Lemma~\ref{lem:decomp}} Fix $d\in \mathbb{Z}$ and for each $k \in \mathbb{N}$ with $-k <d$, define
\[ f_k = \sum_{\substack{ Q \in \mathcal{D} \\ 2^{-k}< \ell(Q) \le 2^d}} \Delta^\mu_Q f + \sum_{\substack{Q \in \mathcal{D} \\ \ell(Q) = 2^d}} E^\mu_Q f.\]
We claim that the sequence $(f_k)$ converges to $f$ pointwise $\mu$-a.e. and in $L^2(\mu).$ Observe that by the Lebesgue differentiation theorem, for $\mu$-a.e. $x$, if $(R_{\ell})$ is a sequence of nested dyadic cubes shrinking to $x$, then
\[  \lim_{\ell \rightarrow \infty} \langle f \rangle^\mu_{R_\ell}  = f(x) \  \ \text{and } \lim_{\ell \rightarrow \infty}  \langle b_Q \rangle^\mu_{R_\ell} =b_Q(x) \qquad \forall Q \in \mathcal{D}.\]
As $\mathcal{S}_{b}$ is $\mu$-sparse, for $\mu$-a.e. $x\in \mathbb{R}^n$, $x$ is in at most finitely many $P \in \mathcal{S}_{b}.$ Then for $\mu$-a.e. $x$, define $P_x$ as follows: if $x$ is in some cube in $\mathcal{S}_{b},$ let $P_x$ denote the smallest cube in $\mathcal{S}_{b}$ containing $x$. Otherwise, let $P_x$ denote any cube containing $x$. Fix any $k$ sufficiently large so that $2^{-k} \le \ell(P_x).$ Then
\begin{equation} \label{eqn:fk}  f_k(x) = \frac{\langle f \rangle_{R_k}^{\mu}}{\langle b_{P_x} \rangle^\mu_{R_k}} b_{P_x}(x),\end{equation}
where $R_k$ is the unique cube containing $x$ with $\ell (R_k) = 2^{-k}$. This, paired with our earlier comments, shows that $(f_k)$ converges $\mu$-a.e.~to $f$. Moreover, \eqref{eqn:fk} implies that $|f_k(x)| \lesssim (M_d^{\mu}f)(x)$ $\mu$-a.e., where $M_d^{\mu}$ is the dyadic maximal function. Then an application of the Dominated convergence theorem gives the $L^2(\mu)$-convergence.

An application of the Cauchy-Schwarz inequality immediately implies that each
 $\| E^\mu_{Q} f \|_{L^2(\mu)} \lesssim \| f 1_{Q} \|_{L^2(\mu)}$ and so to prove the square function estimate, we need only show
$ \sum_{Q \in \mathcal{D}}  \| \Delta^\mu_Q f \|^2_{L^2(\mu)} \lesssim \| f\|^2_{L^2(\mu)}$. To obtain this, we consider
\[
\sum_{Q \in \mathcal{D}}  \| \Delta^\mu_Q f \|^2_{L^2(\mu)}  = \sum_{Q \in \mathcal{D}}  \int_{Q} \Big | \frac{\langle f \rangle_{Q}^{\mu}}{\langle b_{Q}\rangle^{\mu}_{Q}} b_Q - \frac{\langle f \rangle_{Q^{(1)}}^{\mu}}{\langle b_{Q^{(1)}}\rangle^{\mu}_{Q^{(1)}}} b_{Q^{(1)}} \Big |^2 d \mu \lesssim S_1 + S_2,\]
where one can insert $\pm \frac{\langle f \rangle_Q^{\mu}}{ \langle b_{Q^{(1)}}\rangle^{\mu}_{Q^{(1)}}} b_{Q^{(1)}}$ into each integral and estimate the resulting values to get
\[
\begin{aligned}
S_1 &=
\sum_{Q \in S_{b}} |\langle f \rangle^{\mu}_Q|^2 \mu(Q) + |\langle f\rangle_Q^{\mu} - \langle f \rangle_{Q^{(1)}}^{\mu} |^2 \mu(Q),  \\
S_2 &=
\sum_{Q \not \in S_{b}}  |\langle f \rangle^{\mu}_Q|^2 |\langle b_{Q^{(1)}} \rangle_{Q}^{\mu} - \langle b_{Q^{(1)}} \rangle^{\mu}_{Q^{(1)}}|^2 \mu(Q) + |\langle f\rangle_Q^{\mu} - \langle f \rangle_{Q^{(1)}}^{\mu} |^2 \mu(Q).
\end{aligned}
\]
The Carleson embedding theorem paired with the fact that $S_{b}$ is $\mu$-sparse  implies that $\sum_{Q \in S_{b}} |\langle f \rangle^{\mu}_Q|^2 \mu(Q) \lesssim \|f \|_{L^2(\mu)}^2.$ Similarly, Theorem \ref{thm:USF} implies that
\[ \sum_{Q\in \mathcal{D}} | \langle f \rangle_Q^{\mu} - \langle f \rangle^\mu_{Q^{(1)}} |^2 \mu(Q) \lesssim \| f\|^2_{L^2(\mu)}.\]
Thus it remains to bound the first term in $S_2$. To do this, we show that the sequence $(\beta_Q)$ defined by
\begin{equation} \label{eqn:beta} \beta_Q =  |\langle b_{Q^{(1)}} \rangle_{Q}^{\mu} - \langle b_{Q^{(1)}}  \rangle^{\mu}_{Q^{(1)}}|^2 \mu(Q), \quad \text{ for } Q \not \in \mathcal{S}_{b},\end{equation}
and $\beta_Q=0$ otherwise, is a $\mu$-Carleson sequence. To that end, fix a cube $R$ and first assume $R \in \mathcal{D}_{b}$. Then for all $Q \subseteq R$ with $Q \not \in \mathcal{S}_{b}$, we have $Q \in \mathcal{D}_{b}$ and $P_Q = P_{Q^{(1)}}$.  Then using Theorem \ref{thm:USF}, we have
\[
\begin{aligned}
\sum_{Q \subseteq R} \beta_Q & \le |\langle b_{R^{(1)}} \rangle_{R}^{\mu} - \langle b_{R^{(1)}}  \rangle^{\mu}_{R^{(1)}}|^2 \mu(R) + \sum_{\substack{Q\subsetneq R, Q \not \in S_{b} \\ P_Q = P_R}} |\langle b_{Q^{(1)}} \rangle_{Q}^{\mu} - \langle b_{Q^{(1)}}  \rangle^{\mu}_{Q^{(1)}}|^2 \mu(Q) \\
& \hspace{.5in} +
\sum_{P \subsetneq R, P \in \mathcal{S}_{b}} \sum_{\substack{Q: Q \not \in S_{b}  \\ P_Q =P }} |\langle b_{Q^{(1)}} \rangle_{Q}^{\mu} - \langle b_{Q^{(1)}}  \rangle^{\mu}_{Q^{(1)}}|^2\mu(Q) \\
 &\lesssim \mu(R) +  \sum_{Q\subsetneq R}  |\langle b_{P_R} \rangle_{Q}^{\mu} - \langle b_{P_R} \rangle^{\mu}_{Q^{(1)}}|^2 \mu(Q) \\
 & \hspace{.5in} +
\sum_{P \subsetneq R, P \in \mathcal{S}_{b}} \sum_{Q \subsetneq P}  |\langle b_P \rangle_{Q}^{\mu} - \langle b_P \rangle^{\mu}_{Q^{(1)}}|^2 \mu(Q) \\
&\lesssim \mu(R)+  \| b_{P_R} 1_R \|_{L^2(\mu)}^2 + \sum_{P \subsetneq R, P \in \mathcal{S}_{b}} \| b_P\|_{L^2(\mu)}^2 \\
&\lesssim \mu(R) + \sum_{P \subsetneq R, P \in \mathcal{S}_{b}} \mu(P) \\
&\lesssim \mu(R),
\end{aligned}
\]
as needed. Similarly, if $R \in \mathcal{C}_{b}$ then we can write
\[
\begin{aligned}
\sum_{Q \subseteq R} \beta_Q & \le  |\langle b_{R^{(1)}} \rangle_{R}^{\mu} - \langle b_{R^{(1)}}  \rangle^{\mu}_{R^{(1)}}|^2 \mu(R) + \sum_{\substack{Q\subsetneq R, Q \in \mathcal{C}_{b}}} |\langle b_{Q^{(1)}} \rangle_{Q}^{\mu} - \langle b_{Q^{(1)}}  \rangle^{\mu}_{Q^{(1)}}|^2 \mu(Q) \\
&\hspace{.5in}+
\sum_{P \subsetneq R, P \in \mathcal{S}_{b}} \sum_{\substack{Q: Q \not \in S_{b}  \\ P_Q =P }} |\langle b_{Q^{(1)}} \rangle_{Q}^{\mu} - \langle b_{Q^{(1)}}  \rangle^{\mu}_{Q^{(1)}}|^2\mu(Q), \end{aligned}
\]
where first and third terms are bounded as before and the second term equals
\[ \sum_{Q\subsetneq R; Q \in \mathcal{C}_{b}}  |\langle b_R \rangle_{Q}^{\mu} - \langle b_R \rangle^{\mu}_{Q^{(1)}}|^2 \mu(Q) \lesssim \| b_R\|_{L^2(\mu)}^2 \lesssim \mu(R).\]
Thus, $(\beta_Q)$ is $\mu$-Carleson, which completes the proof of estimate (i).  \\

To prove the dual square function estimate  $ \sum_{Q \in \mathcal{D}}  \| (\Delta^\mu_Q)^* f \|^2_{L^2(\mu)} \lesssim \| f\|^2_{L^2(\mu)}$, recall \eqref{eqn:dual}. Then we have
\[
\begin{aligned}
\sum_{Q \in \mathcal{D}} \| (\Delta^\mu_Q)^* f  \|^2_{L^2(\mu)} &= \sum_{Q \in \mathcal{D}} \sum_{Q' \in \ch Q} \left| \frac{ \langle b_{Q'} f \rangle^{\mu}_{Q'}}{\langle b_{Q'} \rangle_{Q'}^\mu} - \frac{\langle b_Q f\rangle_Q^{\mu}}{ \langle b_Q\rangle_Q^{\mu}} \right|^2 \mu(Q') \\
&=
\sum_{Q \in \mathcal{D}} \left| \frac{ \langle b_{Q} f \rangle^{\mu}_Q}{\langle b_{Q} \rangle^\mu_{Q}} - \frac{\langle b_{Q^{(1)}} f\rangle_{Q^{(1)}}^{\mu}}{ \langle b_{Q^{(1)}}\rangle_{Q^{(1)}}^{\mu}} \right|^2 \mu(Q).
\end{aligned}
\]
By inserting $\pm \langle b_{Q^{(1)}} f\rangle_{Q}^{\mu}/\langle b_{Q^{(1)}}\rangle_{Q^{(1)}}^{\mu}$, it is easy to see that this sum is bounded by $S_3 + S_4 +S_5$, where
\[
\begin{aligned}
S_3 &= \sum_{Q \in \mathcal{S}_{b}}  \left(  |\langle b_{Q} f \rangle_Q^\mu|^2 + | \langle b_{Q^{(1)}}f \rangle_Q^{\mu}|^2\right)\mu(Q)  \\
S_4 & = \sum_{Q \not \in \mathcal{S}_{b} }  |\langle b_Q f \rangle_Q^\mu|^2 \left | \langle b_Q\rangle_Q^{\mu} - \langle b_{Q^{(1)}} \rangle_{Q^{(1)}}^{\mu}\right|^2 \mu(Q) \\
S_5 &= \sum_{Q\in \mathcal{D}} \left | \langle b_{Q^{(1)}}  f \rangle_{Q^{(1)}}^{\mu} - \langle b_{Q^{(1)}} f \rangle_Q^{\mu} \right|^2  \mu(Q).
\end{aligned}
\]
Clearly, $S_3 \lesssim \| f\|^2_{L^2(\mu)}$  because $\mathcal{S}_{b}$ is $\mu$-sparse and $  |\langle b_{Q} f \rangle_Q^\mu|^2, | \langle b_{Q^{(1)}}f \rangle_Q^{\mu}|^2 \lesssim \langle |f| \rangle_Q^\mu$. Similarly, $S_4$ is bounded because the sequence $(\beta_Q)$ defined in \eqref{eqn:beta} is $\mu$-Carleson. Thus, we need only consider $S_5.$
Observe that we can decompose $S_5$ as
\[
\sum_{Q: Q^{(1)} \in \mathcal{C}_{b}}  \left | \langle b_{Q^{(1)}}  f \rangle_{Q^{(1)}}^{\mu} - \langle b_{Q^{(1)}} f \rangle_Q^{\mu} \right|^2  \mu(Q) +  \sum_{Q: Q^{(1)} \in \mathcal{D}_{b}}  \left | \langle b_{Q^{(1)}}  f \rangle_{Q^{(1)}}^{\mu} - \langle b_{Q^{(1)}} f \rangle_Q^{\mu} \right|^2  \mu(Q) := S_5^1 + S_5^2.\]
Then Theorem \ref{thm:USF} paired with the properties of $\mathcal{C}_{b}$ give
\[ S_5^1  =\sum_{Q: Q^{(1)} \in \mathcal{C}_{b}}  \left | \langle b  f \rangle_{Q^{(1)}}^{\mu} - \langle b  f \rangle_Q^{\mu} \right|^2  \mu(Q) \lesssim \|b f\|^2_{L^2(\mu)} \lesssim \|f \|^2_{L^2(\mu)},\]
where $b$ is defined in \eqref{eqn:bcb}. To consider $S_5^2$, first write
\begin{equation} \label{eqn:iv2} S_5^2 = \sum_{P \in \mathcal{S}_{b}}  \sum_{Q: P_{Q^{(1)}}=P}  \left | \langle b_P  f \rangle_{Q^{(1)}}^{\mu} - \langle b_P f \rangle_Q^{\mu} \right|^2  \mu(Q).
\end{equation}
For each $P \in \mathcal{S}_{b}$, let $\mathcal{D}_P$ denote the set of maximal $S \in \mathcal{S}_{b}$ so that $S \subsetneq P$. If $Q$ satisfies $P_{Q^{(1)}}=P$,  then for $J = Q, Q^{(1)}$ we can write
\[ \langle b_P  f \rangle_{J}^{\mu}  =\left \langle  1_{(P \setminus \cup_{S \in \mathcal{D}_P} S)} b_P  f \right \rangle_{J}^{\mu}  + \left\langle \sum_{S \in \mathcal{D}_P} 1_S \langle b_P  f \rangle_S^{\mu} \right\rangle_J^\mu. \]
This uses the fact that the $S\in \mathcal{D}_P$ are disjoint and if $S \cap J \ne \emptyset$, then since $S \subsetneq P$ and $P_{Q^{(1)}} = P$, we must have $S \subseteq J.$
Substituting that into $\eqref{eqn:iv2}$ for $J=Q, Q^{(1)}$ and using Theorem \ref{thm:USF} gives
\[
\begin{aligned}S_5^2 &\lesssim \sum_{P \in \mathcal{S}_{b}} \left \|   1_{(P \setminus \cup_{S \in \mathcal{D}_P} S)} b_P  f\right \|^2_{L^2(\mu)} + \sum_{P \in \mathcal{S}_{b}}  \left \|  \sum_{S \in \mathcal{D}_P} 1_S \langle b_P  f \rangle_S^{\mu}  \right \|^2_{L^2(\mu)} \\
& \lesssim \sum_{P \in \mathcal{S}_{b}} \left \|   1_{(P \setminus \cup_{S \in \mathcal{D}_P} S)}  f\right \|^2_{L^2(\mu)} + \sum_{P \in \mathcal{S}_{b}}  \sum_{S \in \mathcal{D}_P} |
\langle b_P  f \rangle_S^{\mu}|^2 \mu(S) \\
& \lesssim \| f\|^2_{L^2(\mu)} + \sum_{S \in \mathcal{S}_{b}} |\langle |f| \rangle_S^\mu|^2 \mu(S) \\
& \lesssim \|f \|^2_{L^2(\mu)},
\end{aligned}
\]
where we use the fact that if $P,R \in \mathcal{S}_{b}$, then the sets $P \setminus \cup_{S \in \mathcal{D}_P} S$ and $R \setminus \cup_{S \in \mathcal{D}_R} S$ are disjoint, and the fact that $\mathcal{S}_{b}$ is $\mu$-sparse. \end{proofof}

\end{document}